\documentclass[11pt, a4paper, reqno]{amsart}
\usepackage{amssymb, array, amsmath, amscd, pdfpages, enumerate, amsthm, fixltx2e, setspace, hyperref,mathtools}

\usepackage[utf8]{inputenc}
\usepackage{amssymb}
\usepackage{bm}
\usepackage{graphicx}

\newcommand*{\R}{{\mathbb{R}}}

\newcommand*{\Abs}[2][default]{\ifthenelse{\equal{#1}{default}}{\left\lvert#2\right\rvert}{\ldelim{#1}{\lvert}#2\rdelim{#1}{\rvert}}}
\newcommand*{\Norm}[2][default]{\ifthenelse{\equal{#1}{default}}{\left\lVert#2\right\rVert}{\ldelim{#1}{\lVert}#2\rdelim{#1}{\rVert}}}

\newcommand*{\Iprod}[3][default]{\ifthenelse{\equal{#1}{default}}{\left\langle#2,#3\right\rangle}{\ldelim{#1}{\langle}#2,#3\rdelim{#1}{\rangle}}}
\newcommand*{\Dualpair}[3][default]{\ifthenelse{\equal{#1}{default}}{\left\langle#2,#3\right\rangle}{\ldelim{#1}{\langle}#2,#3\rdelim{#1}{\rangle}}}

\newcommand*{\ddb}[2][1]{\ifthenelse{\equal{#1}{1}}{\frac{d}{d#2}}{\frac{d^{#1}}{d#2^{#1}}}}
\newcommand*{\pd}[3][1]{\ifthenelse{\equal{#1}{1}}{\frac{\partial{#2}}{\partial{#3}}}{\frac{\partial^{#1}{#2}}{\partial#3^{#1}}}}

\newcommand*{\prk}{P_{\scriptscriptstyle K}}
\newcommand*{\prks}{P_{\scriptscriptstyle K^*}}

\newcommand*\lenv{{\hbox{\raisebox{-.15ex}{\rotatebox[origin=c]{50}{$\smallsmile$}}\kern-8.65pt\rotatebox[origin=c]{-25}{$\smallsetminus$}}}}
\newcommand*\uenv{{\hbox{\raisebox{-.0ex}{\rotatebox[origin=c]{-45}{$\smallfrown$}}\kern-5.4pt\raisebox{.2ex}{\rotatebox[origin=c]{-5}{\scriptsize\slash}}}}\,\kern+1.5pt}

\newcommand{\sleqk}{\preccurlyeq_{\scriptscriptstyle K}}
\newcommand{\sgeqk}{\succcurlyeq_{\scriptscriptstyle K}}
\newcommand{\leqk}{\leq_{\scriptstyle *}}
\newcommand{\geqk}{\geq_{\scriptstyle *}}
\DeclareMathOperator{\mmin}{Min}
\DeclareMathOperator{\mmax}{Max}
\newcommand{\mini}[1]{\mmin [\kern1pt #1 \kern1pt ]}
\newcommand{\maxi}[1]{\mmax [\kern1pt #1 \kern1pt ]}

\newcommand{\sleq}{\preccurlyeq}
\newcommand{\sgeq}{\succcurlyeq}

\newcommand*{\lupp}[1]{\kern0pt^{\kern0pt u \kern0pt}\kern1pt#1}
\newcommand*{\llow}[1]{\kern0pt^{\kern0pt l \kern0pt}\kern1pt#1}
\newcommand*{\rupp}[1]{#1\kern1pt^{\kern0pt u \kern0pt}\kern0pt}
\newcommand*{\rlow}[1]{#1\kern1pt^{\kern0pt l \kern0pt}\kern0pt}
\newcommand*{\ul}[1]{\kern0pt^{\kern0pt u \kern0pt}\kern0pt#1\kern1pt^{\kern0pt l \kern0pt}\kern0pt}
\newcommand*{\lu}[1]{\kern0pt^{\kern0pt l \kern0pt}\kern0pt#1\kern0pt^{\kern0pt u \kern0pt}\kern0pt}
\newcommand*{\bproofname}{Proof}
\newenvironment{bproof}[1][\bproofname]{\begin{proof}[#1]}{\end{proof}}

\newtheorem{thm}{Theorem}[section]
\newtheorem{prop}[thm]{Proposition}

\theoremstyle{definition}
\newtheorem{defn}[thm]{Definition}

\newtheorem{remark}[thm]{Remark}
\newtheorem{example}[thm]{Example}

\numberwithin{equation}{section}

\begin{document}

\title[Mixed lattice structures and cone projections]{Mixed lattice structures and cone projections}

\thispagestyle{plain}

\author{Jani Jokela}
\address[J. Jokela]{Mathematics, Faculty of Information Technology and Communication Sciences, 
Tampere University, PO. Box 692, 33101 Tampere, Finland}
\email{jani.jokela@tuni.fi}
%\address[J. Jokela]{}
%\email{jani.m.jokela@gmail.com}

\begin{abstract}
Problems related to projections on closed convex cones are frequently encountered in optimization theory and related fields. To study these problems, various unifying ideas have been introduced, including  
asymmetric vector-valued norms and certain generalized lattice-like operations. 
We propose a new perspective on these studies by describing how the problem of cone projection can be formulated using an order-theoretic formalism developed in this paper. 
The underlying mathematical structure is a partially ordered vector space that generalizes the notion of a vector lattice by using two partial orderings and having certain lattice-type properties with respect to these orderings. In this note we introduce a generalization of these so-called mixed lattice spaces, and we show how such structures arise quite naturally in some of the applications mentioned above.
\end{abstract}


\maketitle

\section{Introduction}
\label{sec:intro}

The problem of cone projection and the related decomposition method with respect to dual cones is of significant importance in convex analysis,  
and it has been studied extensively since the pioneering work of J. Moreau \cite{moreau}. In this paper we present an order-theoretic setting in which the basic projection problem can be stated and studied. 
The theoretical framework on which this note is based    
is called a \emph{mixed lattice space} which is a partially ordered vector space (or more generally, a group) with two partial orderings and certain ''mixed'' lattice-type properties with respect to these two orders. In a mixed lattice space, the vector lattice operations of supremum and infimum are replaced by asymmetric mixed upper and lower envelopes, which are formed with respect to the two partial orders. A mixed lattice space is a generalization of a vector lattice in the sense that if the two partial orders are identical, then the mixed lattice space reduces to a vector lattice. However, the two partial orders of a mixed lattice space do not need to be lattice orderings. A systematic development of this theory is quite recent \cite{jj1,jj2,jj3}, although some earlier studies on mixed lattice groups date back to 1990's \cite{eri1,eri}. The theory of mixed lattice spaces and groups is heavily based on an earlier theory of mixed lattice semigroups, which was developed by M. Arsove and H. Leutwiler in the 1970's for the purposes of axiomatic potential theory (see \cite{ars} and references therein).

In Section \ref{sec:sec2} we give a brief overview of the terminology and the elementary properties of mixed lattice spaces. In Section~\ref{sec:sec3} we present a  
generalization of mixed lattice structures by relaxing some of the assumptions in the definitions, resulting in a significant gain in generality. 
This modification gives the structure more flexibility in terms of applications to cone projections, and this is what our main results are concerned with. In Sections \ref{sec:sec4} and \ref{sec:sec5} we formulate certain problems in convex optimization using this generalized order-theoretic framework. Here we observe a close connection between the mixed lattice theory and the theory of cone projections and %isotone cone 
retractions. 

In metric cone projections there are two cones involved (a cone and its dual), and we show that this setting has a generalized mixed lattice structure. We then look at some other fundamental facts and results concerning cone projections, and formulate them in the language of mixed lattice theory. We also observe that the so-called \emph{lattice like operations} introduced in \cite{nemeth2013} are in fact a special case of the generalized mixed lattice operations. 

As another application, in Section \ref{sec:sec5} we consider the connections of mixed lattice spaces to the concept of \emph{isotone cone retraction}, introduced by S. N\'emeth in \cite{nemeth2011}, 
and the related notion of an \emph{asymmetric cone norm}, recently studied in \cite{nemeth2020}. Here we show that the mixed envelopes in a mixed lattice space can be used to construct asymmetric cone norms, much in the same way as the positive part mapping in a vector lattice defines an asymmetric cone norm, but our setting is more general.

\section{Mixed lattice spaces}
\label{sec:sec2}

In this section we collect the essential terminology, definitions and basic results. For a more detailed presentation with proofs, we refer to \cite{eri,jj1,jj2}. %\cite{jj1} and \cite{jj2}.

Let $V$ be a partially ordered real vector space %(or an additive group) 
with two partial orderings $\leq$ and $\sleq$ (i.e. we assume that both orderings are compatible with the linear structure %(or additive group structure) 
of $V$). The partial order $\leq$ is called the \emph{initial order} and $\sleq$ is called the \emph{specific order}. 
With these two partial orders $\leq$ and $\sleq$ we define the \emph{mixed upper and lower envelopes} 
\begin{equation}\label{upperenv}
u\uenv v\,=\,\min \,\{\,w\in V: \; w\sgeq u \; \textrm{ and } \; w\geq v \,\}
\end{equation}
and
\begin{equation}\label{lowerenv}
u\lenv v \,=\,\max \,\{\,w\in V: \; w\sleq u \; \textrm{ and } \; w\leq v \,\}, %\quad \textrm{and} \quad
\end{equation}
respectively, where the minimum and maximum (whenever they exist) are taken with respect to the initial order $\leq$. These definitions were introduced by Arsove and Leutwiler in \cite{ars3}. 
We observe that these operations are not commutative, i.e. $x\uenv y$ and $y\uenv x$ are not necessarily equal. 
We recall that a subset $K$ of a vector space is called a \emph{cone} if \, (i)\;$t K \subseteq K$ for all $t \geq 0$, \, (ii)\;\;$K+K \subseteq K$ and \, (iii) \,$K\cap (-K)=\{0\}$. Although the above definition of a cone is rather standard in the theory of partially ordered spaces, we should point out that there are also different definitions in use, and sometimes a subset satisfying the above conditions is called a \emph{convex pointed cone}.

%Def. M.L.Space

\begin{defn}\label{lml}
Let $V$ be a partially ordered real vector space %(or an additive group) 
with respect to two partial orders $\leq$ and $\sleq$, and let $V_p$ and $V_{sp}$ be the corresponding positive cones, respectively. Then $V=(V,\leq,\sleq)$ is called a \emph{mixed lattice vector space} %(or a \emph{mixed lattice group}, respectively) 
if the following conditions hold:
\begin{enumerate}[(1)]
\item
The elements $x\lenv y$ and $x\uenv y$
exist in $V$ for all $x,y\in V$,
\item
$V_{sp}\subseteq V_p$,
\item
%$V_{sp}$ is a \emph{mixed lattice cone}, that is 
The elements $x\uenv y$ and $x\lenv y$ are in $V_{sp}$ whenever $x,y\in V_{sp}$.
\end{enumerate}
\end{defn}

\begin{remark}
The condition $(2)$ can alternatively be stated as follows: if $x,y\in V$ then $x\sleq y$ implies $x\leq y$. 
A cone satistying the condition $(3)$ is a called a \emph{mixed lattice cone}. We also note that the cone $V_{sp}$ does not need to be a generating cone, 
whereas the cone $V_p$ is always generating, that is $V=V_p - V_p$. This is a consequence of Theorem \ref{absval}(b) given below. Various examples of mixed lattice spaces are given in \cite[Section 2]{jj1}.
\end{remark}

The definition of the mixed envelopes implies that in a mixed lattice space $V$ the inequalities 
$x\lenv y \sleq x\sleq x\uenv y$ and $x\lenv y \leq y\leq x\uenv y$ 
hold for all $x,y\in V$. 
For proofs and further discussion on the properties given in the next theorem, see \cite{eri,jj1,jj2}. 

%MUUTA LAUSEEKSI:
\begin{thm}\label{basic1}
Let $V$ be a mixed lattice space. The mixed envelopes have the following properties for all $x,y,z,u,v\in V$ and $a\in \R$.
\begin{enumerate}[(a)]
\item
$x\uenv y \, + \, y \lenv x \, = \, x+y$  %\\ %\quad \textrm{for all} \; x,y\in G.\\
\item
$z \,+ \,x\uenv y \, = \, (x+z)\uenv(y+z)$ \, and \, $z \,+\, x\lenv y \, = \, (x+z)\lenv(y+z)$ %\\
%\item
%$z \,+\, x\lenv y \, = \, (x+z)\lenv(y+z)$ %\\
\item
$x\uenv y \, = \, -(-x \lenv -y)$ %\\
\item
$x\sleq u \quad \textrm{and} \quad y\leq v \; \implies \; x\uenv y \leq u\uenv v \quad \textrm{and} \quad x
\lenv y \leq  u\lenv v$ %\\
%\item
%$x\sleq y \, \implies \, x\leq y$  %\\
\item
$x\leq y \, \iff \, y\uenv x = y \, \iff \, x\lenv y = x$ %\\
\item
$x\sleq y \, \iff \, x\uenv y = y \, \iff \, y\lenv x = x$ %\\
%\item
%$x\sleq u \quad \textrm{and} \quad y\leq v \; \implies \; x\uenv y \leq u\uenv v \quad \textrm{and} \quad x
%\lenv y \leq  u\lenv v$ %\\
%\item
%$x\sleq z \quad \textrm{and} \quad y\leq z \; \implies \; x\uenv y \leq z \quad \textrm{and} \quad x\lenv y \leq z$ % \\
\item
$(ax)\lenv(ay) = a(x\lenv y) \quad \textrm{and} \quad (ax)\uenv(ay) = a(x\uenv y) \quad \textrm{for all } \; a\geq 0$ %\\
\item
$(ax)\lenv(ay) = a(x\uenv y) \quad \textrm{and} \quad (ax)\uenv(ay) = a(x\lenv y) \quad \textrm{for all }\; a< 0$ %\\
\item
$x\sleq y  \; \implies \; z\uenv x \sleq z\uenv y \quad \textrm{and} \quad z\lenv y \sleq  z\lenv y$ %\\
\item
$u\sleq x\sleq z \quad \textrm{and} \quad u\sleq y\sleq z \; \implies \; x\uenv y \sleq z \quad \textrm{and} \quad u \sleq x\lenv y$%\\
%\item
%$z\sleq x \quad \textrm{and} \quad z\sleq y \; \implies \; z \sleq x\lenv y$  %\\%\bigskip\\
%\end{array} 
%$$
\end{enumerate}
\end{thm}

The upper and lower parts of an element were introduced in \cite[Definition 3.1]{jj1}. Their roles are similar to those of the positive and negative parts of an element in a vector lattice.   
 
\begin{defn}\label{upperparts}
Let $V$ be a mixed lattice space and $x\in V$. The elements $\lupp{x}=x\uenv 0$ and $\llow{x}=(-x)\uenv 0$ are called the \emph{upper part} and \emph{lower part} of $x$, respectively. Similarly, the elements $\rupp{x}=0\uenv x$ and $\rlow{x}=0\uenv (-x)$ are called \emph{specific upper part} and \emph{specific lower part} of $x$, respectively. 
\end{defn}

From the above definitions we observe that for the specific upper and lower parts we have $\rupp{x}\sgeq 0$ and $\rlow{x}\sgeq 0$, and for the upper and lower parts $\lupp{x}\geq 0$ and $\llow{x}\geq 0$. 
%The elements $\rupp{x},\rlow{x}\in V_{sp}$ and $\lupp{x},\llow{x}\in V_p$ can be thought of as ''projections'' of $x$ onto the positive cones $V_{sp}$ and $V_p$, respectively.
%
%
The upper and lower parts 
have several important basic properties, which were proved in \cite[Section 3]{jj1}. Some of these properties are given in the next theorem.

\begin{thm}\label{absval} %[{\cite{jj1}} Theorem ?]
Let $V$ be a mixed lattice space and $x\in V$. Then we have
\begin{enumerate}[(a)]
\item
\; $\lupp{x}\,=\,\llow{(-x)}$ \quad \textrm{and} \quad $\rupp{x}\,=\,\rlow{(-x)}$. 
\item
\; $x \, = \, \rupp{x} \, -\,\llow{x} \, = \, \lupp{x} \, - \, \rlow{x}$.
\item
\;
$\lupp{(x+y)}\leq \lupp{x}+\lupp{y}$ \quad and \quad $\rlow{(x+y)}\leq \rlow{x}+\rlow{y}$. 
\item
\;
$\rupp{(x+y)}\leq \rupp{x}+\rupp{y}$ \quad and \quad $\llow{(x+y)}\leq \llow{x}+\llow{y}$. 
\item
\; $\rupp{x}\lenv \llow{x} \, = \, 0\, = \, \rlow{x}\lenv \lupp{x}$.
\item
\; $x\sgeq 0$ \, if and only if \, 
$x=\lupp{x}=\rupp{x}$\, and \,$\llow{x}=\rlow{x}=0$.
\item
\; $x\geq 0$ \, if and only if \, $x=\lupp{x}$\, and \,$\rlow{x}=0$.
\end{enumerate}
\end{thm}

\section{The generalized mixed lattice structure}
\label{sec:sec3}

In a mixed lattice space $V$ 
the existence of the mixed envelopes places rather strict restrictions on the cones $V_p$ and $V_{sp}$ which limits the range of possible applications. 
In this section we present a generalization of a mixed lattice space to overcome this limitation. Our main motivation for the generalization is that when studying cone projections in optimization problems, the upper and lower parts of elements do not typically exist. 

Let $(V,\leq,\sleq)$ be a partially ordered vector space with two partial orderings, as in the preceding section. We introduce the following set notation. 
$$
[x\uenv y]=\{w: w\sgeq x \textrm{ and } w\geq y\} \quad \textrm{and} \quad 
[x\lenv y]=\{w: w\sleq x \textrm{ and } w\leq y\}.
$$
Let $E\subset V$. An element $x\in E$ is called a \emph{minimal element} of $E$ if  
$y\in E$ and $y\leq x$ implies $y=x$. 
A dual notion of a \emph{maximal element} is defined similarly. 
The set of minimal elements of the set $[x\uenv y]$ will be denoted by $\mini{x\uenv y}$, and the set of maximal elements of the set $[x\lenv y]$ will be denoted by $\maxi{x\lenv y}$.

We now define a generalization of the mixed lattice structure in which the elements $x\uenv y$ and $x\lenv y$ are replaced by set-valued mappings $(x,y)\mapsto \mini{x\uenv y}$ and $(x,y)\mapsto \maxi{x\lenv y}$. This provides a considerable increase in generality, at the expense of losing some good properties of mixed lattice spaces and becoming somewhat more difficult to work with.

%mixed quasi-lattice structure:

\begin{defn}\label{gmls}
Let $V$ be a partially ordered real vector space %(or an additive group) 
with respect to two partial orders $\leq$ and $\sleq$, and let $V_p$ and $V_{sp}$ be the corresponding positive cones, respectively. Then $V=(V,\leq,\sleq)$ is called a \emph{generalized mixed lattice structure} if $V_{sp}\subseteq V_p$ and the sets $\mini{x\uenv y}$ and $\maxi{x\lenv y}$ 
are non-empty for all $x,y\in V$. 
\end{defn}

The assumption that the sets $\mini{x\uenv y}$ and $\maxi{x\lenv y}$ are non-empty implies that the cone $V_p$ is generating. 
Clearly, if the sets $\mini{x\uenv y}$ and $\maxi{x\lenv y}$ contain only one element for every $x,y\in V$ then these elements are equal to $x\uenv y$ and $x\lenv y$, respectively, and the generalized mixed lattice structure then reduces to an ordinary mixed lattice space.

In the following we introduce selected fundamental properties of generalised mixed lattice structures to facilitate the study of cone projections in Section~\ref{sec:sec4}. 
The next theorem gives some of the basic properties of generalized mixed lattice structures corresponding to the properties of the mixed envelopes listed in Theorem \ref{basic1}$(a)-(h)$.

\begin{thm}\label{main}
Let $V$ be a generalized mixed lattice structure. The following hold for all $x,y\in V$.
\begin{enumerate}[(a)]
\item
$\mini{x\uenv y}=-\maxi{(-x)\lenv (-y)}$
\item
If $u\in \mini{x\uenv y}$ then there exists an element $w\in \maxi{y\lenv x}$ such that $x+y=u+w$. Similarly, for any $w\in \maxi{y\lenv x}$ there exists an element $u\in \mini{x\uenv y}$ such that $x+y=u+w$. Hence, we have $x+y\in \mini{x\uenv y}+\maxi{y\lenv x}$.
\item
$\mini{(z+x)\uenv (z+y)}=z+\mini{x\uenv y}$ \,for all $z\in V$.
\item 
$\maxi{(z+x)\lenv (z+y)}=z+\maxi{x\lenv y}$ \,for all $z\in V$.
\item
If $x\sleq w$ and $y\leq z$ then for every $u\in \maxi{x\lenv y}$ there exists an element $v\in \maxi{w\lenv z}$ such that $u\leq v$. Similarly, for every $v\in \mini{w\uenv z}$ there exists an element $u\in \mini{x\uenv y}$ such that $u\leq v$.
\item
$x\leq y \; \iff \; \maxi{x\lenv y}=\{x\} \; \iff \; \mini{y\uenv x}=\{y\}$.
\item
$x\sleq y \; \iff \; \maxi{y\lenv x}=\{x\} \; \iff \; \mini{x\uenv y}=\{y\}$.
\item
For all $a\in \R$, $a\geq 0$ we have $\maxi{(ax)\lenv (ay)}=a \maxi{x\lenv y}$ and $\mini{(ax)\uenv (ay)}=a \mini{x\uenv y}$.
\item
For all $a\in \R$, $a< 0$ we have $\maxi{(ax)\lenv (ay)}=a \mini{x\uenv y}$ and $\mini{(ax)\uenv (ay)}=a \maxi{x\lenv y}$.
\end{enumerate}
\end{thm}

\begin{bproof}
\emph{(a)}\, If $u\in \mini{x\uenv y}$ then $u\sgeq x$ and $u\geq y$, so $-u\sleq -x$ and $-u\leq -y$, and thus $-u\in [(-x)\lenv (-y)]$. If $w\in [(-x)\lenv (-y)]$ and $w\geq -u$ then by a similar argument $-w\in [x\uenv y]$. But $-w\leq u$ and $u$ is a minimal element of $[x\uenv y]$, so we must have $w=u$. Hence, $u\in -\maxi{(-x)\lenv (-y)}$. The reverse inclusion is proved similarly, so $\mini{x\uenv y}=-\maxi{(-x)\lenv (-y)}$.

\emph{(b)}\, If $u\in \mini{x\uenv y}$ then $u\sgeq x$ and $u\geq y$. This implies that $x+y\sleq u+y$ and $x+y\leq x+u$, and so $x+y-u\in [y\lenv x]$. Suppose that $w\in [y\lenv x]$ and $w\geq x+y-u$. Then $w\sleq y$ and $w\leq x$, and it follows that $u\geq x+y-w\sgeq x$ and $u\geq x+y-w\geq y$. Hence $u\geq x+y-w\in [x\uenv y]$. But $u$ is a minimal element of the set $[x\uenv y]$, so we must have $u=x+y-w$ and thus $w=x+y-u\in \maxi{y\lenv x}$ and $u+w=x+y$. The dual statement is proved similarly.

\emph{(c)}\, If $v\in \mini{x\uenv y}$ then $v \sgeq x$ and $v \geq y$. Consequently, $v+z \sgeq x+z$ and $v+z \geq y+z$ for all $z\in V$, and so $v+z\in [(x+z)\uenv (y+z)]$. If $w\in [(x+z)\uenv (y+z)]$ and $w\leq v+z$ then $v\geq w-z\sgeq x+z-z=x$ and $v\geq w-z\geq y+z-z=y$. But then $v\geq w-z\in [x\uenv y]$, and since $v$ is a minimal element, it follows that $v=w-z$, or $w=v+z$. Hence, $z+v\in \mini{(z+x)\uenv (z+y)}$ for all $z\in V$.
 
For the converse, if $w\in \mini{(z+x)\uenv (z+y)}$ then $w\sgeq z+x$ and $w\geq z+y$. Thus $w-z\sgeq x$ and $w-z\geq y$, so $w-z\in [x\uenv y]$. Again, if $v\in [x\uenv y]$ and $v\leq w-z$ then $w\geq v+z\sgeq x+z$ and $w\geq v+z\geq y+z$, and it follows that $w\geq v+z\in [x\uenv y]$. Since $w$ is minimal, we have $w=v+z$, or $v=w-z$. Hence, $w-z\in \mini{x\uenv y}$. This shows that $w=z+(w-z)\in z+\mini{x\uenv y}$, proving the equality of the two sets.

\emph{(d)} is similar to \emph{(c)}.

\emph{(e)}\, If $x\sleq w$, $y\leq z$ and $u\in \maxi{x\lenv y}$ then $u\sleq x\sleq w$ and $u\leq y\leq z$. Hence, $u\in [w\lenv z]$ so there is some $v\in \maxi{w\lenv z}$ such that $u\leq v$. The second statement is proved similarly.

\emph{(f)}\, Assume $x\leq y$. Since $y\sleq y$, it follows from $(e)$ that for every $z\in \mini{y\uenv y}=\{y\}$ there exists some $w\in \mini{y\uenv x}$ such that $w\leq z=y$. But $w\sgeq y$ which implies $w\geq y$, and so $w=y$. This shows that $\mini{y\uenv x}=\{y\}$. Furthermore, by $(b)$ there exists $v\in \maxi{x\lenv y}$ such that $w+v=y+v=y+x$, so $v=x$. Hence $\maxi{x\lenv y}=x$. 
Conversely, if $\maxi{x\lenv y}=\{x\}$, or if $\mini{y\uenv x}=\{y\}$, then $x\leq y$.

\emph{(g)} is similar to \emph{(f)}.

\emph{(h)}\, The case $a=0$ is trivial. If $a> 0$ and $z\in\mini{(ax)\uenv (ay)}$ then $z\sleq ax$ and $z\leq ay$. Hence $\frac{z}{a}\sleq x$ and $\frac{z}{a}\leq y$, and so $z=a\frac{z}{a}\in a \mini{x\uenv y}$. The reverse inclusion is straightforward, and the other equality is similar.

\emph{(i)} follows from \emph{(a)} and \emph{(h)}.
\end{bproof}

The converse of the property \emph{(e)} in the preceding theorem does not hold in general, that is, if $x\sleq w$ and $y\leq z$ then there may exist elements $u\in \mini{x\uenv y}$ such that $u\leq v$ does not hold for any $v\in \mini{w\uenv z}$. This is illustrated in the following example.

\begin{example}\label{esim1}
Let $V=\R^3$ and consider the ''pyramid cone''
$$
C=\{a(1,1,1)+b(-1,1,1)+c(-1,-1,1)+d(1,-1,1) :a,b,c,d\geq 0\}.
$$
Let $V_p=V_{sp}=C$. Then $(V,\leq,\sleq)$ is a generalized mixed lattice. Next, let $x=(2,2,0)$ and $y=(-2,2,0)$. Then $\mini{x\uenv y}=\{(0,t,2):0\leq t\leq 4\}$. If $w=(0,4,2)$ and $z=(0,2,2)$ then $\mini{w\uenv z}=\{(s,3,3):-1\leq s\leq 1\}$, $w\sgeq x$ and $z\geq y$. Now if $u=(0,0,2)$ then $u\in \mini{x\uenv y}$ %$u\sgeq x$ and $u\geq y$ 
but there is no such element $v\in \mini{w\uenv z}$ that $v\geq u$.
\end{example}

In fact, the property %$(e)$ 
discussed above characterizes mixed lattice spaces.

\begin{thm}
Let $(V,\leq,\sleq)$ be a generalized mixed lattice structure with elements $x,y,z,w\in V$ such that $x\sleq w$ and $y\leq z$. The following statements are equivalent. 
\begin{enumerate}[(a)]
\item
$V$ is a mixed lattice space.
\item
For every $u\in \mini{x\uenv y}$ there exists an element $v\in \mini{w\uenv z}$ such that $u\leq v$.
\item
For every $v\in \maxi{w\lenv z}$ there exists an element $u\in \maxi{x\lenv y}$ such that $u\leq v$.
\end{enumerate}
\end{thm}

\begin{bproof}
The implication $(a) \implies (b)$ follows by the definition of mixed lattice space. Assume that $(b)$ holds and let $x,y\in V$. Suppose there are elements $a,u\in \mini{x\uenv y}$ and $a\neq u$. Then $a,u\sgeq x$ and $a,u\geq y$ so it follows by hypothesis that, in particular, there exists an element $v\in \mini{a\uenv a}=\{a\}$ such that $a=v\geq u$. But $a$ and $u$ are both minimal elements, so they are comparable only if they are equal, that is $u=a$. This is a contradiction, so the set $\mini{x\uenv y}$ contains only one element $x\uenv y$. This shows that $V$ is a mixed lattice space. 
The equivalence $(a) \iff (c)$ is proved similarly. 
\end{bproof}

The next result shows that in a generalized mixed lattice structure every element can be written as a difference of a positive part and a negative part, but the representation is not unique. 
The following is thus a generalized version of \cite[Theorem 3.6]{jj1} and parts (b) and (e) of Theorem \ref{absval}.

\begin{thm}\label{representation}
Let $V$ be a generalized mixed lattice structure and $x\in V$. 
\begin{enumerate}[(a)]
\item
For any $u\in \mini{x\uenv 0}$ there exist an element $w\in\mini{0\uenv (-x)}$ 
such that $x=u-w$. Then $V_{p}\cap\maxi{w\lenv u}=\{0\}$ %$0\in\maxi{w\lenv u}$ 
and $u+w\in\mini{u\uenv w}$. 
Conversely, if $x=u-w$ and  $0\in\maxi{w\lenv u}$ 
then $u\in\mini{x\uenv 0}$ and $w\in\mini{0\uenv (-x)}$. 
\item
For any $u\in \mini{0\uenv x}$ there exist an element $w\in \mini{(-x)\uenv 0}$ such that $x=u-w$. Then $V_{p}\cap\maxi{u\lenv w}=\{0\}$ and $u+w\in \mini{w\uenv u}$. 
Conversely, if $x=u-w$ and $0\in\maxi{u\lenv w}$ then $u\in\mini{0\uenv x}$ and $w\in\mini{(-x)\uenv 0}$.
\end{enumerate}
\end{thm}

\begin{bproof}
\emph{(a)} and \emph{(b)} are similar so we only prove \emph{(a)}.  
The first part follows immediately from Theorem \ref{main}(b). Indeed, if $x\in V$ and $u\in \mini{x\uenv 0}$ then there exists $v\in \maxi{0\lenv x}$ such that $x=u+v$. If we put $w=-v$ then $w\in \mini{0\uenv (-x)}$ and $x=u-w$. 

Next we note that $0\in \maxi{0\lenv 0}=\{0\}$, and since $0\sleq w$ and $0\leq u$, it follows by Theorem \ref{main}(e) that there exists an element $y\in \maxi{w\lenv u}$ such that $y\geq 0$. Now $y\in \maxi{w\lenv u}=\maxi{w\lenv (x+w)}=w+\maxi{0\lenv x}$ (by Theorem \ref{main}(d)) so there exists $z\in \maxi{0\lenv x}$ such that $y=w+z$. But then $z=y-w\geq -w$, and since $-w$ and $z$ are both maximal elements of the set $[0\lenv x]$, they are comparable only if $z=-w$, and so $y=0$. This shows that $\maxi{w\lenv u}\cap V_{p}=\{0\}$ %$0\in \maxi{w\lenv u}$ 
and so by Theorem \ref{main}(b) $u+w\in\mini{u\uenv w}$. 

Conversely, if $x=u-w$ and $0\in\maxi{w\lenv u}$ then by Theorem \ref{main}(c) we have $u+w\in\mini{u\uenv w}=\mini{(x+w)\uenv w}=w+\mini{x\uenv 0}$. Hence, $u\in\mini{x\uenv 0}$. On the other hand, we have $u=x+w\in\mini{x\uenv 0}=x+\mini{0\uenv (-x)}$ and this shows that $w\in \mini{0\uenv (-x)}$, finishing the proof. 
\end{bproof}

By the preceding theorem, the set $\mini{0\uenv x}$ can be called the \emph{set of specific upper parts} of $x$, and $\mini{(-x)\uenv 0}$ the \emph{set of lower parts} of $x$. For any $x\in V$ we can choose an upper part $u\in\mini{0\uenv x}$, and there always exists a corresponding lower part $v\in\mini{(-x)\uenv 0}$ such that $x=u-v$. Similar remarks apply to the sets $\mini{x\uenv 0}$ and $\mini{0\uenv (-x)}$, called the \emph{set of upper parts} of $x$ and the \emph{set of specific lower parts} of $x$, respectively.

%simplify notation: by theorem \ref{presentation} we can simplify the notation and denote by $[\rupp{x}]=\mini{0\uenv x}$ jne, called the sets of upper parts and lower parts. For every $x$ we can choose any upper part $u\in[\rupp{x}]$ and the corresponding lower part $v\in[\llow{x}]$ such that $x=u-v$.

%tähän: dual cones

We introduce some additional terminology for the next section.  
The set 
$$
V_{sp}^*=\{y\in V: V_p \cap \maxi{x\lenv y}\neq \emptyset \textrm{ for all } x\in V_{sp} \}
$$ 
is called the \emph{right dual of} $V_{sp}$, and the set
$$
^* V_p=\{x\in V: V_p \cap \maxi{x\lenv y}\neq \emptyset\textrm{ for all } y\in V_p \}
$$ 
is called the \emph{left dual of} $V_{p}$. 

In fact, the cones $V_{sp}$ and $V_p$ are duals of each other. 

\begin{prop}\label{dualcones2}
$V_{sp}^*=V_p$ and $^*V_p=V_{sp}$. Hence, $V_{sp}= ^{*}\!\!(V_{sp}^*)$ and $(^*V_p)^*=V_{p}$.
\end{prop}

\begin{bproof}
Let $y\in V_{sp}^*$. Then for any $x\in V_{sp}$ there exist $w\in\maxi{x\lenv y}$ and $0\leq w\leq y$. Thus $y\in V_p$. On the other hand, if $y\in V_p$ then by Theorem \ref{main}(e) for any $x\in V_{sp}$ there exist $w\in\maxi{x\lenv y}$ such that $0\leq w$. Hence $y\in V_{sp}^*$, proving that $V_{sp}^*=V_p$. 
Next, if $x\in V_{sp}$ and $y\geq 0$ then, again by Theorem \ref{main}(e) there exist $w\in\maxi{x\lenv y}$ such that $v\geq 0$, and so $x\in ^{*}\!\!V_p$. Let $x\in ^{*}\!\!V_p$. Then for any $y\in V_{p}$ there exist $w\in\maxi{x\lenv y}$ such that $w\geq 0$. In particular, if $y=0$ then there exist $v\in\maxi{x\lenv 0}$ such that $v\geq 0$. But then $0\leq v\leq 0$, so $v=0$ and $v\sleq x$ implies that $x\in V_{sp}$. This shows that $^*V_p=V_{sp}$. 
\end{bproof}

\section{Mixed lattice structure in the problem of cone projection}
\label{sec:sec4}

The main application of the results of this paper are given in this section. 
The aim is to 
show how the problem of cone projection 
can be stated in a purely order-theoretic form in the framework of generalized mixed lattice structure,  
thus providing a new perspective on such problems. 
We assume the knowledge of basic notions and terminology of convex optimization. For these we refer to \cite{boyd} and \cite{fca}.

%\subsection{Cone projections} \hfill\\

Let $K$ be a closed and convex pointed cone in $\R^n$ with the dual cone $K^*=\{y:\langle x , y \rangle \geq 0 \textrm{ for all } x\in K\}$. Let $\sleqk$ be the partial ordering induced by the cone $K$ 
and let $\leqk$ be the partial ordering given by the dual cone $K^*$. 
We may assume here that $K\subseteq K^*$. Indeed, if this is not the case, then (since $(K^*)^*=K$) we can simply exchange the roles of $K$ and $K^*$ in the following discussion by setting $C=K^*$ and $C^*=K$. Then $C\subseteq C^*$, and $\preccurlyeq_{\scriptscriptstyle C}$ is defined as the order induced by $C$ and $\leqk$ is the order given by $C^*$. 

Let $\prk :\R^n \to K$ be the projection mapping that gives the unique point $\prk  x$ on $K$ nearest to $x$. That is,
$$
\prk  x\in K \quad \text{and} \quad ||x-\prk  x||=\inf \{||x-y||:y\in K\}.
$$
This nearest point $\prk  x$ has the characterization (\cite[Theorem 3.1.1]{fca})
\begin{equation}\label{nearpoint}
\prk  x\in K \quad \text{and} \quad \langle \prk  x -x , \prk  x -y \rangle \leq 0 \textrm{ for all } y\in K.
\end{equation}
The projection mapping also has the translation property
\begin{equation}\label{translate}
P_{x+K}y = x+\prk  (y-x) \quad \textrm{ for all } x,y\in \R^n.
\end{equation}
The %projection 
mapping $\prk $ is called $K$-\emph{isotone} if $x\sleqk y$ implies $\prk  x\sleqk \prk  y$.

Clearly, the projection $\prk  x$ satifies $\prk  x\sgeqk 0$ and $\prk  x\geqk x$. 
%eka suoraan proj.määr. toka karakterisoinnista 4.2 sillä jos y=0 niin sisätulo on 0 (koska proj ja proj-x ovat ortog.) joten duaalikartion määr.mukaan proj >0 
Now we can show that $V=(\R^n , \leqk, \sleqk)$ is a generalized mixed lattice structure in the sense of Definition \ref{gmls}, and the orthogonal projection $\prk  x$ is in fact a minimal element satisfying the inequalities $\prk  x\sgeqk 0$ and $\prk  x\geqk x$. %, that is, $\prk  x\in \mini{0\uenv x}$.

\begin{thm}\label{lemma1}
Let $K$ be a closed and convex cone in $\R^n$ and $K^*$ its dual cone such that $K\subseteq K^*$, and let $\sleqk$ and $\leqk$ be the partial orderings determined by the cones $K$ and $K^*$, respectively.  
Then $V=(\R^n , \leqk,\sleqk)$ is a generalized mixed lattice structure and for every $x\in \R^n$ the projection element $\prk x$ satisfies $\prk  x\in \mini{0\uenv x}$. %Similarly, the projection $\prks x$ satisfies $\prks  x\in \mini{x\uenv 0}$.
\end{thm} 

\begin{bproof}
Let $x\in \R^n$. Since $\prk x$ is the orthogonal projection of $x$ on $K$, we have $\langle \prk  x -x , \prk  x \rangle = 0$, and from this we get $\langle \prk  x , \prk  x \rangle = \langle x , \prk  x  \rangle$. As noted above, the element $\prk  x$ satifies $\prk  x\sgeqk 0$ and $\prk  x\geqk x$. Suppose there is some other element $w$ such that  $w\sgeqk 0$, $w\geqk x$ and $w\leqk \prk  x$. Then $\prk  x\in K$ and $\prk  x -w\geqk 0$, or  $\prk  x -w\in K^*$, so by the definition of $K^*$ we have $\langle \prk  x -w, \prk  x  \rangle \geq 0$, so $\langle \prk  x , \prk  x  \rangle\geq \langle w , \prk  x  \rangle$. On the other hand, $w-x\in K^*$, and so $\langle w- x , \prk  x  \rangle\geq 0$, which gives $\langle w , \prk  x  \rangle\geq \langle  x , \prk  x  \rangle$. Hence, we have 
$$
\langle \prk  x , \prk  x \rangle = \langle x , \prk  x  \rangle \leq \langle w , \prk  x \rangle \leq \langle \prk  x , \prk  x  \rangle.
$$ 
Thus, 
$\langle \prk  x , \prk  x \rangle = \langle x , \prk  x  \rangle = \langle w , \prk  x \rangle$, and it follows that 
$$
\langle \prk  x , \prk  x \rangle - \langle w , \prk  x \rangle = \langle \prk  x -w , \prk  x \rangle=0.
$$ 

Now, for every $y\in K$ we have 
$$
\langle \prk  x -w , \prk  x - y \rangle = \langle \prk  x -w , \prk  x  \rangle - \langle \prk  x -w , y \rangle,
$$
where the first term is zero, as we have shown, and for the second term we have $\langle \prk  x -w , y \rangle \geq 0$ by the definition of the dual cone. Hence, we have $\langle \prk  x -w , \prk  x - y \rangle\leq 0$ for all $y\in K$. By characterization \eqref{nearpoint} this means that $\prk  x$ is the unique point on $K$ nearest to $w$. But $w\in K$, so we must have $w=\prk  x$. This shows that $\prk  x\in \mini{0\uenv x}$, and since this holds for any $x\in V$ it follows by Theorem \ref{main} that the sets $\mini{x\uenv y}$ and $\maxi{x\lenv y}$ are non-empty for all $x,y\in V$, and so $V$ is a generalized mixed lattice structure. %Proof of the second statement regarding the element $\prks x$ is similar.
\end{bproof}

Now Theorem \ref{lemma1} allows us to translate the projection problem to the mixed lattice setting.  
Let us choose an element $\rupp{x}\in \mini{0\uenv x}$ by the criterion of shortest distance, that is, $\rupp{x}=\prk  x$. Then, if we denote by $\llow{x}$ the corresponding element in $\mini{(-x)\uenv 0}$ (see Theorem \ref{representation}), we have the unique representation $x=\rupp{x}-\llow{x}$ for every $x\in V$, and this is the most natural representation in the present setting.

After fixing the ''representatives'' $\rupp{x}$ and $\llow{x}$ of each $x$ in this way, we can now simplify (or rather abuse) the notation and write $\rupp{x}=0\uenv x$ 
and $\llow{x}=(-x)\uenv 0=-(x\lenv 0)$ (we observe here that the element $\llow{x}$ gives the projection of $-x$ on $K^*$).  
It now follows from the results of Section~\ref{sec:sec3} that, in essence, our generalized mixed lattice structure behaves much like an ordinary mixed lattice space, and (with some care) we can apply the rules $(a)-(h)$ of Theorem \ref{basic1}. 
Hence, by Theorem \ref{main}(c) equation \eqref{translate} becomes 
\begin{equation}\label{translate2}
P_{x+ \scriptscriptstyle K}y = x+\prk  (y-x)=x+\rupp{(y-x)}=x+ 0\uenv(y-x)=x\uenv y. 
\end{equation} 
In other words, $x\uenv y$ is the point on the cone $x+K$ that is nearest to the point $y$, or equivalently, the point on the cone $y+K^*$ that is nearest to the point $x$. 
In a similar manner, the lower envelope $x\lenv y$ is associated with the projections on the cones $x-K$ and $y-K^*$, i.e. $x\lenv y = P_{x- \scriptscriptstyle K}y = P_{y- \scriptscriptstyle K^*}x$.

We observe that $K$ and $K^*$ are mutually dual cones also in the order-theoretic sense of Proposition \ref{dualcones2}. We next show that orthogonality in the usual sense implies the order-theoretic version of the orthogonality condition.

\begin{prop}\label{ortho}
If $x\in K$ and $y\in K^*$ are elements such that $\langle x,y\rangle=0$ then $0\in\maxi{x\lenv y}$.
\end{prop}

\begin{bproof}
Let $z\in K$. Then by the definition of $K^*$ we have $\langle z,y\rangle\geq 0$, and since $\langle x,y\rangle=0$ we obtain
$\langle x-z,y\rangle=\langle x,y\rangle - \langle z,y\rangle=- \langle z,y\rangle \leq 0$. 
This holds for all $z\in K$, so if we define $u=x-y$ then $y=x-u$ and the above inequality becomes $\langle x-u,x-z\rangle \leq 0$ for all $z\in K$. By the characterization \eqref{nearpoint} this means that $x=\prk u$. Then by Theorem \ref{lemma1} we have $x\in \mini{0\uenv u}$, so by Theorem \ref{main}(c) we get 
$y=x-u\in \mini{0\uenv u}-u=\mini{(-u)\uenv 0}$. 
Hence, $0\in\maxi{x\lenv y}$, by Theorem \ref{representation}.
\end{bproof}

A fundamental tool in the study of cone projections is the following classical theorem of Moreau \cite{moreau}.

\begin{thm}[Moreau]\label{moreau1} 
Let $K$ be a closed convex cone in $\R^n$ and $K^*$ its dual cone. 
Every $x\in \R^n$ can be written as  
$x=\prk x- \prks (-x)$ %(or $x=\lupp{x}-\rlow{x}$) 
where $\langle \prk x , \prks (-x) \rangle=0$. Moreover, $\prk x=0$ holds if and only if $x\in -K^*$. 
\end{thm}

Using the notation introduced above, we get the following special case of Theorem \ref{representation}, 
which can be viewed as the order-theoretic version of Moreau's theorem.

\begin{thm}\label{moreau2}
Let $V=(\R^n , \leqk,\sleqk)$ be the generalized mixed lattice structure on $\R^n$, where $\sleqk$ is the partial order defined by a closed convex cone $K$ and $\leqk$ is the partial order defined by the dual cone $K^*$ such that $K\subseteq K^*$. Then every $x\in V$ can be written as $x=\rupp{x}-\llow{x}$ %(or $x=\lupp{x}-\rlow{x}$) 
where $\rupp{x}\sgeqk 0$, $\llow{x}\geqk 0$ and $\rupp{x}\lenv \llow{x}=0$. 
Moreover, $\rupp{x}=0$ if and only if $x\leqk 0$. %$x\in -K^*$.
\end{thm}

\begin{bproof}
By Theorem \ref{lemma1}, for any $x\in V$ we have $\prk x\in\mini{0\uenv x}$, and
as noted above, if we put $\rupp{x}=\prk x$ then $\rupp{x}\sgeqk 0$ and by Theorem \ref{representation}, if $\llow{x}$ is the corresponding element in $\mini{(-x)\uenv 0}$, then $\llow{x}=\prks (-x)\geqk 0$ and we have the unique representation $x=\rupp{x}-\llow{x}$.  
Since $\langle \rupp{x},\llow{x}\rangle=0$, it follows by Proposition \ref{ortho} that $0\in \maxi{\rupp{x}\lenv \llow{x}}$. 
Also, by Theorem \ref{main} we have $\maxi{\rupp{x}\lenv \llow{x}}=-\mini{(-\rupp{x})\uenv (-\llow{x})}$, so from equation \eqref{translate2} we get $P_{-\rupp{x}+ \scriptscriptstyle K} (-\llow{x})=-(\rupp{x}\lenv \llow{x})$. On the other hand, by \eqref{translate} we obtain
$$
P_{-\rupp{x}+ \scriptscriptstyle K} (-\llow{x})=-\rupp{x}+\prk (-\llow{x}+\rupp{x})=-\rupp{x}+\prk  x=-\rupp{x}+\rupp{x}=0,
$$
and this gives justification for writing $\rupp{x}\lenv \llow{x}=0$. (Again, this  
just amounts to the fact that we have chosen the %choosing an appropriate 
''representative'' from the set $\maxi{\rupp{x}\lenv \llow{x}}$ by the criterion of shortest distance, which is consistent with our earlier choice of $\rupp{x}$, that is, if $\rupp{x}$ and $\llow{x}$ are chosen as above then the element in the set $\maxi{\rupp{x}\lenv \llow{x}}$ corresponding to this choice is $0$.)

Finally, by Theorem \ref{main}(f) we have 
$x\leqk 0$ if and only if $\mini{0\uenv x}=\{0\}$, and so $\prk x=\rupp{x}=0$ holds if and only if $x\leqk 0$ (that is, $x\in -K^*$).
\end{bproof}

It should be stressed again that the notation we use here is not entirely correct because our structure is not a mixed lattice space in the sense of Definition \ref{lml}, but rather a generalized structure as described in the preceding section. Although the element $\prk  x$ is the minimum in terms of distance to the cone $K$, it is not necessarily the order-theoretic minimum in the sense of \eqref{upperenv} (this kind of situation occurs in Example \ref{esim1}). Because of this, some properties of mixed lattice space do not hold in the present situation, and some care should be taken when manipulating expressions that contain the mixed envelopes. For instance, the inequalities in Theorem \ref{basic1}(i) do not necessarily hold, which means that the projection mapping is not isotone, in general. The conditions for the isotonicity of the cone projection have been extensively studied (see \cite{nemeth2013} and the references therein). However, the following discussion gives further justification for the use of this notation. 

The authors in \cite{nemeth2013} introduced what they called the ''lattice-like operations'' for studying questions related to cone projections. These operations are a generalization of similar operations that were introduced in \cite{gowda} for the special case of self-dual cones (i.e. $K=K^*$). The lattice-like operations are defined by
$$
x\sqcup y =P_{x+ \scriptscriptstyle K}y, \quad x\sqcap y=P_{x- \scriptscriptstyle K}y, \quad x\sqcup_* y =P_{x+ \scriptscriptstyle K^*}y, \quad x\sqcap_* y=P_{x- \scriptscriptstyle K^*}y. 
$$ 
However, we can observe that these are exactly the mixed lattice operations: 
$$
x\sqcup y=x\uenv y, \quad x\sqcup_* y=y\uenv x, \quad x\sqcap y=x\lenv y \quad \text{and} \quad x\sqcap_* y=y\lenv x.
$$ 
Indeed, by definition we have 
$$
x\sqcup y =P_{x+ \scriptscriptstyle K}y =x+\prk (y-x)=x+0\uenv (y-x)=x\uenv y,
$$ 
and the other identities can be verified similarly. 

This shows that the lattice-like operations in \cite{nemeth2013} are a special case of the generalized mixed lattice operations, and consequently, most of the properties of the lattice-like operations (\cite[Lemma 2 and Lemma 3]{nemeth2013}) 
are identical to the  
properties of the mixed envelopes listed in Theorem \ref{basic1} $(a)-(h)$.  
%
%For example, 
%the orthogonality property $\langle x-x\sqcap y, x\sqcup_* y -x \rangle =0$ follows, since 
%$$x-x\sqcap y=x-x\lenv y=0\uenv (x-y) =\rupp{(x-y)}
%$$ 
%and 
%$$
%x\sqcup_* y -x=y\uenv x-x=(y-x)\uenv 0=\llow{(x-y)}.
%$$ 
%The orthogonality condition $\rupp{(x-y)}\lenv \llow{(x-y)} =0$ then follows from ??.
%
%
It appears that the authors in \cite{nemeth2013} were not aware of the mixed lattice theory, and they discovered these operations independently in this particular special setting. The above discussion places the lattice-like operations in their proper order-theoretic context.

%asymmetric vector norms

\section{On cone retractions and asymmetric cone norms} %\hfill\\
\label{sec:sec5}

Another concept closely related to the cone projections is the notion of a cone retraction, which was introduced by S. N\'emeth in \cite{nemeth2011}, and since then they have been an active area of research. 
The study of such mappings can be further generalized by introducing the notion of an asymmetric cone norm.

\begin{defn}\label{asymm_cone_norm}
Let $K\subset X$ be a cone in a topological vector space $X$ and let $\leq$ be the associated order relation. A continuous mapping $Q:X\to K$ is an \emph{asymmetric cone norm} if 
\begin{enumerate}[(1)]
\item
$Q$ is a retraction onto $K$, that is, $Q(x)=x$ for all $x\in K$ and $Q(X)=K$
\item
$Q(tx)=tQx$ for all $t\in \R_+$ and $x\in X$
\item
$Q(x+y)\leq Qx+Qy$ for all $x,y\in X$
\item
If $Qx=0$ and $Q(-x)=0$ then $x=0$.
\end{enumerate}
Moreover, we say that $Q$ is a \emph{proper asymmetric cone norm} if $Q(I-Q)=0$. Here $I$ denotes the identity operator on $X$.
\end{defn}

The problem of existence of asymmetric cone norms has been studied in \cite{nemeth2020}. 
It is well known that asymmetric cone norms exist for positive cones in topological vector lattices, since the positive part mapping $x\mapsto x^+$ has all the properties given in Definition \ref{asymm_cone_norm}.  Mixed lattice spaces 
provide a more general setting in which asymmetric cone norms exist, and they indeed  
arise very naturally from the mixed lattice order structure. 
The following result is from \cite[Theorem 5.2]{jj3}, and it is a rather immediate consequence of the properties of the upper parts and the mixed envelopes, given in Theorems \ref{basic1} and \ref{absval}.

\begin{thm}[{\cite{jj3}} Theorem 5.2]\label{thm2}
Let $V=(V,\leq,\sleq)$ be a mixed lattice space with a topology such that the mappings 
$x\mapsto \rupp{x}$ and $x\mapsto \lupp{x}$ are continuous. 
Then the following hold.
\begin{enumerate}[(a)]
\item
Let $V_{p}$ be the positive cone associated with the partial order $\leq$. Then the mapping $Q:V\to V_{p}$ given by $Q(x)=\lupp{x}$ is a proper asymmetric cone norm on $V$. 
\item
Let $V_{sp}$ be the positive cone associated with the partial order $\sleq$. Then the mapping $Q:V\to V_{sp}$ given by $Q(x)=\rupp{x}$ is a proper asymmetric cone norm on $V$ which is isotone with respect to both partial orderings.
\end{enumerate} 
\end{thm}

%\begin{proof}
%\emph{(a)}\, It follows at once from the properties of $\lupp{x}$ %given in Theorem \ref{absval} and \ref{basic1}(g) 
%that the mapping $Q(x)=\lupp{x}$ has the properties listed in Definition \ref{asymm_cone_norm}. Property $(1)$ follows by Theorem \ref{absval}(g), and $(2)$ follows by Theorem \ref{basic1}(g). Property $(3)$ is an immediate consequence of Theorem \ref{absval}(c), while $(4)$ follows from Theorem \ref{absval}(a) and (g). 
%To check that $Q$ is proper, we note that 
%$$
%Q(I-Q)(x)=Q(x-\lupp{x})=Q(-\rlow{x})=\lupp{(-\rlow{x})} =\llow{(\rlow{x})} =0,
%$$ 
%by Theorem \ref{absval} (since $\rlow{x}\sgeq 0$). 
%The proof of \emph{(b)} is similar. The isotonicity in \emph{(b)} follows by Theorem \ref{basic1}(d) and (g). 
%\end{proof}

\begin{remark} 
Note that $\llow{x}$ and $\rlow{x}$ are also proper asymmetric cone norms that map $V$ onto $-V_p$ and $-V_{sp}$, respectively. In the preceding theorem we simply assumed the continuity of the mappings $x\mapsto \rupp{x}$ and $x\mapsto \lupp{x}$. 
This is because the continuity of the mixed lattice operations depends on the topology of the space, and we have not considered topologies in this paper. The topological theory of mixed lattice spaces is studied in \cite{jj3}, but this topic is too broad to be discussed here.
\end{remark}

Incidentally, the result in part (b) of the preceding theorem has some other interesting implications.  
For example, the following result was given in \cite{nemeth2011}.

\begin{thm}[{\cite{nemeth2011}} Theorem 2]\label{thm3}
Let $H$ be a Hilbert space and $K\subset H$ a pointed closed convex generating normal cone. If there exists a continuous isotone retraction $f:H\to K$ such that $(f-I)(H)\cap (I-f)(H)=\{0\}$, then $K$ is a lattice cone.
\end{thm}

It is interesting to note that in the mixed lattice theory a conclusion similar to the above theorem follows from Theorem \ref{thm2}(b), at least in the finite-dimensional case. Indeed, it was proved in \cite[Theorem 4.11]{jj3} that if $V$ is a finite dimensional normed mixed lattice space such that the cone $V_{sp}$ is generating, and the mixed lattice operations are continuous then $V_{sp}$ is a lattice cone. Combining this with Theorem \ref{thm2}(b) yields a result analogous to the above theorem in finite-dimensional mixed lattice spaces.

On the other hand, Theorem \ref{thm2}(a) is related to a more general situation in which the cone retraction is not necessarily isotone. 
Suppose we have a cone $K$ in some topological vector space $X$, and we ask if there exists an asymmetric cone norm $Q$ on $X$ such that $Q(X)=K$ and $Q$ satisfies the conditions in Definition \ref{asymm_cone_norm}. 
Theorem \ref{thm2}(a) shows that a sufficient condition for the existence of such mapping is that there exists a mixed lattice order structure on $X$ such that $K$ is a positive cone for the partial order $\leq$. This is obviously another nontrivial problem in general, 
but we can obtain some interesting results using the mixed lattice techniques. For example, 
it is possible to show that given any closed convex cone in $\R^n$ we can always turn $\R^n$ into a mixed lattice space $(\R^n,\leq,\sleq)$ such that the given cone will be the positive cone associated with the partial order $\leq$, and so by Proposition \ref{thm2} there always exists a continuous asymmetric cone norm associated with the given cone in $\R^n$. The latter part (that is, the existence of a proper asymmetric cone norm) has been proved recently in \cite[Theorem 2]{nemeth2020}. A different proof based on the mixed lattice theory can be found in \cite[Theorem 5.3]{jj3}. These observations further illuminate  
the connection of the mixed lattice theory to the problems related to retraction cones.

\section{Conclusions}

In this paper we presented a new approach to the problem of cone projection based on an ordered algebraic structure called a mixed lattice space.  
We first introduced a generalization of the notion of mixed lattice space, and we showed that many of the basic properties of mixed lattice spaces can be extended to the generalized mixed lattice structure. The  motivation for this generalization is that it can be applied in a broader range of situations. As our main application, we showed how the mixed lattice structure arises quite naturally in the study of cone projections. We demonstrated how the problem of cone projection can be formulated in the mixed lattice setting, and we also observed that the related notion of lattice-like operations can be viewed as a special case of the generalized mixed lattice operations. We also discussed the more general concepts of cone retractions and asymmetric cone norms, and we showed that asymmetric cone norms can be constructed in mixed lattice spaces using the upper part mapping, which is a generalization of the positive part mapping on vector lattices.

%\vspace{0.7cm}
%\noindent
%\small{\textbf{Data availability} \\ 
%Data sharing not applicable to this article as no datasets were generated or analysed during the current study.}
%\noindent
%\small{\textbf{Funding and competing interests} \\ %This study was funded by Tampere University. 
%The author has no funding or competing interests to declare that are relevant to the content of this article.}

%%------------------------------------------------------------------------------------------------

\bibliographystyle{plain}
%\bibliography{reference}

\end{document}